\documentclass[12pt]{article}
\usepackage{amsfonts}

\newcommand{\supp}{{\rm supp}}

\title{Weak Continuity of Dynamical Systems for the KdV and mKdV Equations}
\author{Shangbin Cui$^\dag$ \ \ and \ \ Carlos E. Kenig$^\ddag$\\ [0.1cm]
 {\small $^\dag$ Department of Mathematics, Sun Yat-Sen University,
  Guangzhou, Guangdong 510275}\\
 {\small People's Republic of China. E-mail: cuisb@yahoo.com.cn}\\ [0.1cm]
 {\small $^\ddag$ Department of Mathematics, University of Chicago,
  Chicago, IL 60637, USA.}\\
 {\small E-mail: cek@math.uchicago.edu}}
\date{}

\begin{document}
\maketitle

\begin{abstract}
In this paper we study weak continuity of the dynamical systems for
the KdV equation in $H^{-3/4}(\mathbb{R})$ and the modified KdV
equation in $H^{1/4}(\mathbb{R})$. This topic should have
significant applications in the study of other properties of these
equations such as finite time blow-up and asymptotic stability and
instability of solitary waves. The spaces considered here are
borderline Sobolev spaces for the corresponding equations from the
viewpoint of the local well-posedness theory. We first use a variant
of the method of \cite{GoM} to prove weak continuity for the mKdV,
and next use a similar result for a mKdV system and the generalized
Miura transform to get weak continuity for the KdV equation.
\end{abstract}

\section{Introduction}

\hskip 2em
The purpose of this paper is to establish weak continuity of the
dynamical system $S(t)$ in the Sobolev space $H^{-3/4}(\mathbb{R})$
of the KdV equation
$$
  \partial_t u+\partial_x^3 u-6u\partial_xu=0,
\eqno{(1.1)}
$$
and the dynamical system $S_1(t)$ in $H^{1/4}(\mathbb{R})$ of the
defocusing mKdV equation
$$
  \partial_t u+\partial_x^3 u-6u^2\partial_xu=0.
\eqno{(1.2)}
$$
Here, by the notion that {\em $S(t)$ is a dynamical system} of the KdV
equation (1.1) in some Sobolev space $H^s(\mathbb{R})$ we mean that
$\{S(t):t\in\mathbb{R}\}$ is a family of bounded and continuous (nonlinear)
operators in $H^s(\mathbb{R})$ satisfying the following three conditions:

\begin{description}
\item[$(a)$] $\{S(t):t\in\mathbb{R}\}$ is a strongly continuous group of
bounded and continuous (nonlinear) operators in the space $H^s(\mathbb{R})$
(so that for any $u_0\in H^s(\mathbb{R})$, we have that $[t\mapsto
S(t)u_0]\in C(\mathbb{R},H^s(\mathbb{R}))$, $S(0)u_0=u_0$, and
$S(t)S(t')u_0=S(t+t')u_0$ for all $t,t'\in\mathbb{R}$. Note that since
$S(t)S(-t)=S(-t)S(t)=id$, we have $S(t)^{-1}=S(-t)$).

\item[$(b)$] For every $u_0\in H^s(\mathbb{R})$, the function
$u$ on $\mathbb{R}^2$ defined by $u(x,t)=[S(t)u_0](x)$ (for $(x,t)\in
\mathbb{R}^2$) is a solution (in a weak sense) of the initial value problem
$$
\left\{
\begin{array}{l}
  \partial_t u+\partial_x^3 u-6u\partial_xu=0, \quad
  x\in\mathbb{R},\;\;\; t\in\mathbb{R},\\
  u(x,0)=u_0(x), \quad x\in\mathbb{R}.
\end{array}
\right.
\eqno{(1.3)}
$$
See Remark 1.1 below for more discussions on this condition.

\item[$(c)$] The mapping $u_0\mapsto [t\mapsto S(t)u_0]$ from
$H^s(\mathbb{R})$ to $C(\mathbb{R},H^s(\mathbb{R}))$ is
bounded and uniformly continuous on bounded sets.
\end{description}
For $S_1(t)$ the definition is similar, except for replacing
(1.3) with the problem
$$
\left\{
\begin{array}{l}
  \partial_t u+\partial_x^3 u-6u^2\partial_x u=0, \quad
  x\in\mathbb{R},\;\;\; t\in\mathbb{R},\\
  u(x,0)=u_0(x), \quad x\in\mathbb{R}.
\end{array}
\right.
\eqno{(1.4)}
$$
\medskip

{\bf Remark 1.1}\ \ We note that the condition $(b)$ implicitly implies
that if $s$ is not sufficiently large, say, $s<0$, then the function
$u(x,t)=[S(t)u_0](x)$ must have certain additional regularity
beyond what a general $C(\mathbb{R},H^s(\mathbb{R}))$ class function
possesses so that the product $6u\partial_xu$ or its alternative form
$3\partial_x(u^2)$ makes sense in the distribution category. In a
concrete construction of a dynamical system $S(t)$ this is usually
achieved by building the solution $u(x,t)$ of (1.3) in a subspace of
$C(\mathbb{R},H^s(\mathbb{R}))$ with functions which have more
regularity so that the product $6u\partial_xu$ is meaningful in the
distribution sense. Choice of such a subspace is diversified and may
be different for different situations, cf., e.g., \cite{B}, \cite{CKSTT},
\cite{Guo}, \cite{KPV1} and \cite{KPV3}.
\medskip

In general, for a given partial differential equation of the evolutionary
type, if we can find a function space $X$ in the space variable and a
family of operators $\{S(t):t\in\mathbb{R}\}$ in $X$ satisfying similar
conditions as $(a)$--$(c)$ above with $H^{-3/4}(\mathbb{R})$ replaced by
$X$ and (1.3) replaced by the corresponding initial value problem for
that equation, then we say that $S(t)$ is a dynamical system of that
equation in $X$, and we also say that the inital value problem for that
equation is {\em globally well-posed} in the space $X$. If $S(t)$ is not
defined for all $t\in\mathbb{R}$, but for each $u_0\in X$ there exists a
corresponding $T>0$ such that $S(t)u_0$ is well-defined for $t\in (-T,T)$
and conditions similar to $(a)$--$(c)$ above with obvious modifications are
satisfied, then we say that equation is {\em locally well-posed} in $X$,
and in this case we call $S(t)$ a {\em local dynamical system}.

Existence of $S(t)$ for KdV and $S_1(t)$ for mKdV in Sobolev spaces has
been a goal of study for many years, cf. \cite{B}, \cite{CKSTT},
\cite{Guo}--\cite{K3}, \cite{KPV1}, \cite{KPV3} and references therein,
culminating in the following result:
\medskip

{\bf Theorem 1.2}\ \ {\em $(i)$ For every $s\geq -3/4$, the KdV equation
$(1.1)$ is globally well-posed in $H^s(\mathbb{R})$. Moreover, for any
$s<-3/4$, the KdV equation $(1.1)$ is not locally well-posed in
$H^s(\mathbb{R})$ in the sense that the solution operator, if it exists, is
not uniformly continuous on bounded sets.

$(ii)$ For every $s\geq 1/4$, the mKdV equation $(1.2)$ is globally
well-posed in $H^s(\mathbb{R})$. Moreover, for any $s<1/4$, the mKdV equation
$(1.1)$ is not locally well-posed in $H^s(\mathbb{R})$ in the sense that the
solution operator, if it exists, is not uniformly continuous on bounded sets.}
\medskip

Thus, for every $s\geq -3/4$ there exists a global dynamical system of the
KdV equation in  $H^s(\mathbb{R})$, and for every $s\geq 1/4$ there exists
a global dynamical system of the mKdV equation in $H^s(\mathbb{R})$. These
dynamical systems are unique. Indeed, we have:
\medskip

{\bf Lemma 1.3}\ \ {\em For every $s\geq -3/4$, the dynamical system of the
KdV equation in $H^s(\mathbb{R})$ is unique, and for every $s\geq 1/4$, the
dynamical system of the mKdV equation in $H^s(\mathbb{R})$ is unique.}
\medskip

{\em Proof}:\ \ By a standard energy estimate argument, we can easily prove
that the solution of the problem (1.3) in $C(\mathbb{R},H^s(\mathbb{R}))$
is unique provided $u_0\in H^s(\mathbb{R})$ and $s>3/2$. In fact, by Zhou
\cite{Zhou} we know that this unconditional uniqueness is actually true for
all $s\geq 0$. This immediately implies that the dynamical system for the KdV
equation in $H^s(\mathbb{R})$ is unique provided $s\geq 0$. Since
$H^\infty(\mathbb{R})=\cap_{s\geq 0} H^s(\mathbb{R})$ is dense in
$H^s(\mathbb{R})$ for any $s\in\mathbb{R}$, by using the condition $(c)$ in
the definition of dynamical systems, we easily see that, for any $-3/4\leq
s<0$, if $S^1(t)$ and $S^2(t)$ are two dynamical systems of the KdV equation
in the space $H^s(\mathbb{R})$, then $S^1(t)u_0=S^2(t)u_0$ for all $u_0\in
H^s(\mathbb{R})$ and $t\in\mathbb{R}$, which is exactly the desired assertion.
The proof for the assertion $(ii)$ is similar, because the energy method
also ensures that the solution of the problem (1.4) is also unique in
$C(\mathbb{R},H^s(\mathbb{R}))$ when $s>3/2$. $\quad\Box$

Thus, for every $s>-3/4$, the dynamical system of the KdV equation in
$H^s(\mathbb{R})$ is the restriction of $S(t)$ in $H^{-3/4}(\mathbb{R})$ to
$H^s(\mathbb{R})$, and for every $s>1/4$, the dynamical system of the mKdV
equation in $H^s(\mathbb{R})$ is the restriction of $S_1(t)$ in
$H^{1/4}(\mathbb{R})$ to $H^s(\mathbb{R})$.

{\bf Remark 1.4}\ \ We note that uniqueness of the dynamical system for an
evolution equation in a Banach space $X$ does not mean uniqueness of the
solution for the initial value problem of that equation in the space
$C(\mathbb{R},X)$, because it does not exclude the possibility of existence
of solutions which cannot be expressed in forms of orbits of the dynamical
system.

By the condition $(c)$ we see that if $S(t)$ is a dynamical system of an
evolution equation in some Sobolev space $H^s(\mathbb{R})$ then the mapping
$u_0\mapsto [t\mapsto S(t)u_0]$ from $H^s(\mathbb{R})$ to $C(\mathbb{R},
H^s(\mathbb{R}))$ is continuous.
A natural question is: {\em Is the dynamical system also weakly continuous?}
If the evolution equation under consideration is linear, then the answer to
this question is trivially positive, because we know that every continuous
linear operator in Banach spaces is also weakly continuous. Since, however,
we are considering nonlinear equations, this question by no means has an
obvious answer. Our motivation of asking this question is inspired by the
important series of works of Martel and Merle \cite{MM1}--\cite{MM4}, where
the authors studied finite time blow-up and asymptotic stability and
instability of solitary waves for the generalized KdV equations. One key step
in their strategy in these works is a reduction to a nonlinear Liouville type
theorem, which was further reduced into a corresponding linear one, involving
the linearized operator around the solitary wave. It is in both these steps
that the weak continuity of the flow map for generalized KdV in suitable
Sobolev spaces plays a central role.

Recently, Kenig and Martel \cite{KM} studied the asymptotic
stability of solitons for the Benjamin-Ono equation in
$H^{1/2}(\mathbb{R})$, following the program initiated by Martel and
Merle for the generalized KdV. Thus, a key step in \cite{KM} is to
establish weak continuity of the dynamical system of the BO in
$H^{1/2}(\mathbb{R})$. The proof is very simple and reduces matters
to the uniform continuity of the dynamical system in spaces of
strictly smaller indices. This reduction relies on the fact that BO
is well-posed in $L^2(\mathbb{R})$ (cf. \cite{IK}). The same method
shows this weak continuity for KdV in $H^s(\mathbb{R})$ with
$s>-3/4$ and mKdV in $H^s(\mathbb{R})$ with $s>1/4$. We can then ask
if the dynamical systems $S(t)$ of KdV in $H^{-3/4}(\mathbb{R})$ and
$S_1(t)$ of the mKdV in $H^{1/4}(\mathbb{R})$ are weakly continuous?
Note that since the KdV and mKdV equations do not have a uniformly
continuous flow map when restricted to bounded sets in the spaces
$H^s(\mathbb{R})$ with $s<-3/4$ and $s<1/4$, respectively, the
approach used in \cite{KM} does not work in these critical cases.

Weak continuity of dynamical systems in critical Sobolev spaces
which are critical from the viewpoint of local well-posedness was
first studied by Goubet and Molinet in the reference \cite{GoM},
where the cubic nonlinear Schr\"{o}dinger equation on the line was
studied. For this equation the global well-posedness in
$L^2(\mathbb{R})$ was established in \cite{Ts}, while in \cite{KPV4}
(focusing case) and \cite{CCT} (defocusing case) it was shown that
the flow map is not uniformly continuous in any Sobolev space of
negative index. Thus, the weak continuity in $L^2(\mathbb{R})$ of
the flow map cannot be treated by the approach reviewed in the above
paragraph. Goubet and Molinet \cite{GoM} affirmatively settled this
problem by taking advantage of the ``local smoothing'' effect
estimates together with a suitable uniqueness result.

We would also like mention two recent interesting preprints by L.
Molinet \cite{M1, M2}, which disprove the weak continuity of the
flow maps in $L^2(\mathbb{T})$ for both the cubic Nonlinear
Schr\"{o}dinger equation and the Benjamin-Ono equation, though we
know that the initial value problems of these equations are globally
well-posed in $L^2(\mathbb{T})$.

In this paper we give a positive answer to the weak continuity
question for KdV and mKdV. More precisely, the main purpose of this
paper is to prove the following results:
\medskip

{\bf Theorem 1.5}\ \ {\em The dynamical system $S(t)$ of the KdV equation
$(1.1)$ in $H^{-3/4}(\mathbb{R})$ is weakly continuous for any fixed $t\in
\mathbb{R}$. In fact, we have the following stronger assertion: Assume that
$u_{0n}\in H^{-3/4}(\mathbb{R})$ $(n=1,2,\cdots)$ and $u_{0n}\to u_0$ weakly
in $H^{-3/4}(\mathbb{R})$ as $n\to\infty$. Let $u_n(x,t)=[S(t)u_{0n}](x)$
$(n=1,2,\cdots)$ and $u(x,t)=[S(t)u_0](x)$. Then for any $T>0$ and any
$\varphi\in H^{-3/4}(\mathbb{R})$ we have
$$
  \lim_{n\to\infty}\sup_{|t|\leq T}
  |(u_n(\cdot,t)-u(\cdot,t),\varphi)_{H^{-3/4}(\mathbb{R})}|=0.
\eqno{(1.5)}
$$
}

{\bf Theorem 1.6}\ \ {\em The dynamical system $S_1(t)$ of the mKdV equation
$(1.2)$ in $H^{1/4}(\mathbb{R})$ is weakly continuous for any fixed $t\in
\mathbb{R}$. In fact, we have also the following stronger assertion: Assume
that $u_{0n}\in H^{1/4}(\mathbb{R})$ $(n=1,2,\cdots)$ and $u_{0n}\to u_0$
weakly in $H^{1/4}(\mathbb{R})$ as $n\to\infty$. Let $u_n(x,t)=[S_1(t)u_{0n}](x)$
$(n=1,2,\cdots)$ and $u(x,t)=[S_1(t)u_0](x)$. Then for any $T>0$ and any
$\varphi\in H^{1/4}(\mathbb{R})$ we have
$$
  \lim_{n\to\infty}\sup_{|t|\leq T}
  |(u_n(\cdot,t)-u(\cdot,t),\varphi)_{H^{1/4}(\mathbb{R})}|=0.
\eqno{(1.6)}
$$
}

As in \cite{GoM} and \cite{CK}, we shall use some compactness
arguments together with suitable uniqueness results to prove the
above results. The proof of Theorem 1.6 is easier than that of
Theorem 1.5. The idea of the proof of Theorem 1.6 (following
\cite{GoM} in a simplified situation) is as follows: If a sequence
of solutions $\{u_n\}$ of the equation (1.2) is bounded in
$C([-T,T],H^{1/4}(\mathbb{R}))$, then $\{\partial_t u_n\}$ is
bounded in $C([-T,T],H^{-11/4}(\mathbb{R}))$, so that $\{u_n\}$ has
a subsequence which is strongly convergent in $L^2([-R,R]\times
[-T,T])$. By this fact and a certain uniqueness result, the desired
conclusion follows. See Section 3 for details of the proof. This
argument clearly does not apply to the equation (1.1) (because here
we deal with Sobolev spaces of negative index). Thus, to prove
Theorem 1.5 we shall appeal to the generalized Miura transform
introduced by Christ, Colliander and Tao in \cite{CCT} to reduce the
problem into the corresponding problem for a mKdV system, for which
the above argument applies. See Section 4 for details.
\medskip

{\bf Acknowledgement}\ \ In the previous version of this manuscript
(posted on arXiv: 0909.0794) the proofs of the above theorems used a
different approach which relies on some smoothing effect estimates
and, therefore, are unnecessarily complicated. The approach used
here was suggested to us by Professor L. Molinet. We are glad to
acknowledge our sincere thanks.

The work of the first author is partially supported by the National
Natural Science Foundation of China under the grant number 10771223
as well as a fund from the Sun Yat-Sen University, and was performed
while visiting the University of Chicago under financial support of
China Scholarship Council. He would like to express his thanks to
the Department of Mathematics of the University of Chicago for its
hospitality during his visit. The second author is supported in part
by NSF grant DMS-0456583.

\section{Review of proofs of Theorem 1.2 and the Miura transform}

\hskip 2em
In order to prove Theorems 1.5 and 1.6, we need to have a basic knowledge
about the proofs of Theorem 1.2 and the Miura transform. In this section we
recall these materials.

Global well-posedness of the KdV initial value problem (1.3) in the Sobolev
space $H^{-3/4}(\mathbb{R})$ was established recently by Guo \cite{Guo}
in the framework of the function space $\bar{F}^s$ ($s\geq -3/4$), which is a
dyadic Bourgain-type space with modifications in the low frequency part of
functions by considering the smoothing effect estimate of the Airy equation.
Similar spaces of this type have previously been used by some other authors,
cf. \cite{IK}, \cite{T2}, \cite{Tar} and references therein. In \cite{Guo}
the author first used a contraction mapping argument in the space
$\bar{F}^{-3/4}$ to get local well-posedness of (1.3) in
$H^{-3/4}(\mathbb{R})$, and next he used the
$I$-operator introduced by Colliander et al in \cite{CKSTT} to establish
almost conservation of a modified energy quantity which ensures that the
local solution can be extended into a global one. The function space
$\bar{F}^s$ ($s\geq -3/4$) is defined as follows. Let $\eta_0:\mathbb{R}\to
[0,1]$ denote an even function supported in $[-8/5,8/5]$ and equal to $1$ in
$[-5/4,5/4]$. For $k\in\mathbb{Z}$, $k\geq 1$, let $\eta_k(\xi)=
\eta_0(2^{-k}\xi)-\eta_0(2^{-k+1}\xi)$. We also denote, for all $k\in\mathbb{Z}$,
$\chi_k(\xi)=\eta_0(2^{-k}\xi)-\eta_0(2^{-k+1}\xi)$. It follows that
$$
  \sum_{k=0}^\infty\eta_k(\xi)=1\quad \mbox{for}\;\;
  \xi\in\mathbb{R},
$$
and
$$
  \sum_{k=-\infty}^\infty\chi_k(\xi)=1\quad \mbox{for}\;\;
  \xi\in\mathbb{R}\backslash\{0\}.
$$
Note that $\supp\chi_k\subseteq [-(8/5)2^k,-(5/8)2^k]\cup
[(5/8)2^k,(8/5)2^k]$ for all $k\in\mathbb{Z}$, and
$\supp\eta_k\subseteq [-(8/5)2^k,-(5/8)2^k]\cup [(5/8)2^k,(8/5)2^k]$
for $k\geq 1$. For $k\in\mathbb{N}$ we denote
$$
  I_k=[-2^{k+1},-2^{k-1}]\cup [2^{k-1},2^{k+1}],
$$
and let $X_k$ be the function space
$$
\begin{array}{rl}
  X_k=&\{f\in L^2(\mathbb{R}\times\mathbb{R}): f\;\;\mbox{supported in}\;\;
  I_k\times\mathbb{R}\;\; \mbox{and}\\ [0.2cm]
  &\|f\|_{X_k}:=\sum_{j=0}^\infty 2^{j/2}\|\eta_j(\tau-\xi^3)
   f(\xi,\tau)\|_{L^2_{\xi,\tau}}<\infty\}.
\end{array}
$$
The function space $\bar{F}^s$ is defined as follows:
$$
  \bar{F}^s\!=\!\{u\in\! S'(\mathbb{R}\!\times\!\mathbb{R}):
  \|u\|_{\bar{F}^s}^2:=\!\sum_{k=1}^\infty\! 2^{2sk}\|\eta_k(\xi)
  \widetilde{u}(\xi,\tau)\|_{X_k}^2\!+\!\|\mathcal{F}[\eta_0(\xi)
  \widetilde{u}(\xi,\tau)]\|_{L_x^2 L_t^\infty}^2<\infty\},
$$
where $\widetilde{u}=\mathcal{F}(u)$ represents Fourier transform of
$u$ (in two variables). It can be easily shown that
$$
  \bar{F}^{-3/4}\subseteq C(\mathbb{R},H^{-3/4}(\mathbb{R}))\cap
  L_x^\infty(\mathbb{R},L_{t,{\rm loc}}^2(\mathbb{R}))
  \subseteq C(\mathbb{R},H^{-3/4}(\mathbb{R}))\cap
  L_{\rm loc}^2(\mathbb{R}^2),
$$
so that for any $u\in\bar{F}^{-3/4}$, $u^2$ makes sense.
We let $\bar{F}^s_T$ be the restriction of $\bar{F}^s$ on $\mathbb{R}\times
[-T,T]$, i.e., $u\in\bar{F}^s_T$ if and only if there exists $w\in\bar{F}^s$
such that $w|_{\mathbb{R}\times [-T,T]}=u$, and the norm
$$
  \|u\|_{\bar{F}^s_T}:=\inf\{\|w\|_{\bar{F}^s}:w\in\bar{F}^s,
  \;w|_{\mathbb{R}\times [-T,T]}=u\}.
$$
By using the method of first establishing a bilinear estimate in the space
$\bar{F}^s$ to get a local solution and next using the $I$-operator to
prove that the norm $\|u(\cdot,t)\|_{H^{-3/4}(\mathbb{R})}$ grows only
polynomially fast so that it cannot blow-up in finite time, Guo \cite{Guo}
proved the following result:
\medskip

{\bf Theorem 2.1}\ \ {\it There exists a bounded and locally Lipschitz
continuous mapping
$\Psi: H^{-3/4}(\mathbb{R})\to C(\mathbb{R},H^{-3/4}(\mathbb{R}))$ such that
$(i)$ for any $u_0\in H^{-3/4}(\mathbb{R})$, $t_0\in\mathbb{R}$ and $T>0$,
the function $u=\Psi(u_0)$ belongs to $\bar{F}^{-3/4}_T$ when restricted to
$\mathbb{R}\times [t_0-T,t_0+T]$, and $(ii)$ $u$ is a solution $($in
distribution sense$)$ of the initial value problem $(1.3)$, and it is the
unique solution of $(1.3)$ satisfying the property ensured by $(i)$.
Moreover, there exists a constant $C>0$ such that
$$
  \|u(\cdot,t)\|_{H^{-3/4}(\mathbb{R})}\leq
  C(1+|t|)\|u_0\|_{H^{-3/4}(\mathbb{R})} \quad
  \mbox{for all}\;\;t\in\mathbb{R}.
\eqno{(2.1)}
$$
}

Combining this result with the global well-posedness of (1.3) in $H^s(\mathbb{R})$
for $s>-3/4$ (cf \cite{CKSTT} and \cite{KPV3}) and the result of \cite{CCT}
which states that the solution operator (if it exists) of (1.3) is not locally
uniformly continuous in $H^s(\mathbb{R})$ for $s<-3/4$, we see that the assertion
$(i)$ of Theorem 1.2 follows.

However, we are unable to directly use Theorem 2.1 to prove Theorem
1.5 by following the approach of \cite{CK} and \cite{GoM}. The
reason is that, though we are able to establish estimates of the
form
$$
  \|D_x^{\theta}u\|_{L_x^\infty L_T^2}\leq C\|u\|_{\bar{F}^{-3/4}_T}
$$
for $0\leq\theta<1/4$, unfortunately we are unable to get an
integral estimate for $D_t^{\theta}u$ by $\|u\|_{\bar{F}^{-3/4}_T}$
even locally in both $x$ and $t$, no matter how small that
$\theta>0$ is, and even worse, the equation (1.1) does not seem to
help for a such estimate (these are crucial techniques used in
\cite{CK} and \cite{GoM}). To prove Theorem 1.5 we shall have to
appeal to the so called {\em generalized Miura transform}, which we
shall recall later in this section.

We now turn our attention to the mKdV equation (1.2). Local well-posedness of
the problem (1.4) in $H^s(\mathbb{R})$ with $s\geq 1/4$ was established by
Kenig, Ponce and Vega in \cite{KPV1}. For our purpose we recall this result
in the critical case $s=1/4$. For any time interval $I=[t_0,t_0+T]$, let $X=
X(\mathbb{R}\times I)$ denote the function space defined as follows: We first
introduce a norm $\|\cdot\|_{X(\mathbb{R}\times I)}$ for measurable
functions $u$ on $\mathbb{R}\times I$:
$$
  \|u\|_{X(\mathbb{R}\times I)}:=\|u\|_{L_{t,I}^\infty H_x^{1/4}}+
  \|u\|_{L_x^4 L_{t,I}^\infty}+\|J_x^{1/4}u\|_{L_x^5 L_{t,I}^{10}}
  +\|\partial_x u\|_{L_x^{20} L_{t,I}^{5/2}}
  +\|\partial_x J_x^{1/4}u\|_{L_x^\infty L_{t,I}^2}.
$$
Next we define
$$
  X(\mathbb{R}\times I)=\{u\in C(I,H^{1/4}(\mathbb{R})):
  \|u\|_{X(\mathbb{R}\times I)}<\infty\}.
$$
In \cite{KPV1}, it was proved that for any $u_0\!\in\! H^{1/4}(\mathbb{R})$
there exists corresponding $T\!=\!T(\|u_0\|_{H^{1/4}(\mathbb{R})})\\ >0$ such that the
problem (1.4) has a unique solution in the space $X(\mathbb{R}\times I)$. Thus,
the dynamical system $S_1(t)$ of the mKdV equation (1.2) in $H^{1/4}(\mathbb{R})$
is well-defined at least locally, and $S_1(t)u_0\in X(\mathbb{R}\times [0,T])$
for any $u_0\in H^{1/4}(\mathbb{R})$, where
$T=T(\|u_0\|_{H^{1/4}(\mathbb{R})})$. To show that $S_1(t)$ is actually
defined for all $t\in\mathbb{R}$ we need the {\em Miura transform} $v=M(u)$,
which is defined by
$$
  v=\partial_x u+u^2.
\eqno{(2.2)}
$$
Indeed, by an argument of Colliander, Keel, Staffilani, Takaoka and Tao
\cite{CKSTT}, the Miura transform can be used to prove global well-posedness
of the mKdV equation in $H^{1/4}(\mathbb{R})$ from that of the KdV equation
in $H^{-3/4}(\mathbb{R})$ ensured by Theorem 2.1. For our purpose we review
a few more details of this argument in the following paragraphs.

We first write:
\medskip

{\bf Lemma 2.2}\ \ {\it For any $s\geq 0$, the Miura transform $M$ is a
bounded and continuous mapping from $H^s(\mathbb{R})$ to $H^{s-1}(\mathbb{R})$,
and it is Lipschitz continuous when restricted to any bounded set in
$H^s(\mathbb{R})$. Moreover, if $s\geq 1/4$ then it is injective, and for any
subset $S$ of $H^s(\mathbb{R})$, if $S$ is bounded in $L^2(\mathbb{R})$ and
$M(S)$ is bounded in $H^{s-1}(\mathbb{R})$, then $S$ is bounded in
$H^s(\mathbb{R})$, or more precisely, there exists constant $C>0$ such that
$$
  \|u\|_{H^s}\leq C(\|M(u)\|_{H^{s-1}}+\|u\|_2^2)
\eqno{(2.3)}
$$
for all $u\in H^s(\mathbb{R})$.}
\medskip

{\em Proof}:\ \ From the proofs of Lemma 9.1 in \cite{CKSTT} and Lemma 9.1 in
\cite{CCT}, we easily see that for any $0\leq s<1$ and $u_1,u_2\in
H^s(\mathbb{R})$, the following inequality holds:
$$
  \|M(u_1)-M(u_2)\|_{H^{s-1}(\mathbb{R})}\leq C(\|u_1\|_{H^s(\mathbb{R})}
  +\|u_2\|_{H^s(\mathbb{R})})\|u_1-u_2\|_{H^s(\mathbb{R})}.
$$
Since $H^s(\mathbb{R})$ is an algebra when $s>1/2$, the above inequality is
trivially true for any $s\geq 1$. Thus, for any $s\geq 0$, $M$ is a bounded and
continuous mapping from $H^s(\mathbb{R})$ to $H^{s-1}(\mathbb{R})$, and it is
Lipschitz continuous when restricted to any bounded set in $H^s(\mathbb{R})$.
The assertion that $M$ is injective when $s\geq 1/4$ follows from Lemma 2.3
below, and the last assertion follows from the proof of Lemma 9.2 in
\cite{CKSTT}. $\quad\Box$
\medskip

{\bf Lemma 2.3}\ \ {\it Let $s\geq 1/4$, $a\in L^4(\mathbb{R})$ and $w\in
H^s(\mathbb{R})$. Assume that
$$
  w'(x)+a(x)w(x)=0 \quad \mbox{for}\;\; x\in\mathbb{R}
$$
$($in distribution sense$)$. Then $w=0$.}
\medskip

{\em Proof}:\ \ We fist note that $w\in H^s(\mathbb{R})$ and $s\geq 1/4$
implies that $w\in L^4(\mathbb{R})$. Thus, $aw\in L^2(\mathbb{R})$. To prove
$w=0$ we only need to show that for any $\varphi\in C_0^\infty(\mathbb{R})$,
$$
  \int_{-\infty}^\infty w(x)\varphi(x)dx=0.
$$
Let $\varphi_0\in C_0^\infty(\mathbb{R})$ be such that $\varphi_0(x)=1$ for
$|x|\leq 1$, $\varphi_0(x)=0$ for $|x|\geq 2$, and $0\leq\varphi_0\leq 1$.
Set $\varphi_n(x)=\varphi_0(x/n)$, $n=1,2,\cdots$. Then $\|\varphi_n\|_\infty
=1$, $n=1,2,\cdots$, and $\|\varphi_n'\|_2\leq Cn^{-1/2}\to 0$ as $n\to\infty$.
Given $\varphi\in C_0^\infty(\mathbb{R})$, let
$$
  \psi(x)=\int_{-\infty}^x\!\!\varphi(y)e^{\int_y^x a(t)dt}dy \quad \mbox{and}
  \quad  \psi_n(x)=\varphi_n(x)\psi(x), \quad n=1,2,\cdots.
$$
We have that $\psi\in L^\infty(\mathbb{R})\cap H_{{\rm loc}}^1(\mathbb{R})$,
$\psi_n\in H_c^1(\mathbb{R})$ (i.e., $\psi_n\in H^1(\mathbb{R})$ and has
compact support), $n=1,2,\cdots$, and
$$
  \psi_n'-a\psi_n=\varphi_n(\psi'-a\psi)+\varphi_n'\psi
  =\varphi_n\varphi+\varphi_n'\psi\to\varphi \quad
  \mbox{in}\;\; L^2(\mathbb{R})
$$
(as $n\to\infty$). Hence
$$
\begin{array}{rcl}
  \displaystyle\int_{-\infty}^\infty w(x)\varphi(x)dx
  &=&\displaystyle\lim_{n\to\infty}\int_{-\infty}^\infty
  w(x)[\psi_n'(x)-a(x)\psi_n(x)]dx
\\ [0.3cm]
  &=&-\displaystyle\lim_{n\to\infty}\int_{-\infty}^\infty
  [w'(x)+a(x)w(x)]\psi_n(x)dx=0.
\end{array}
$$
$\Box$

As well-known, if $u$ is a solution of (1.2) then its Miura transform $v=M(u)$
is a solution of the KdV equation (1.1). Thus, by using a similar argument as
in the proof of Lemma 1.3 we conclude that
$$
  M(S_1(t)u_0)=S(t)(M(u_0))
\eqno{(2.4)}
$$
for any $u_0\in H^{1/4}(\mathbb{R})$ and any $t\in\mathbb{R}$ such that
$S_1(t)u_0$ is well-defined (note that $S(t)u_0$ is well-defined for all
$t\in\mathbb{R}$ whereas so far we only know that $S_1(t)u_0$ is well-defined
for small $|t|$). Using this
relation, the $L^2$ conservation law for the mKdV equation, the growth
estimate (2.1) for solutions of the KdV equation, and Lemma 2.2, we can easily
infer that for any $u_0\in H^{1/4}(\mathbb{R})$, the solution $u=S_1(t)u_0$
of the initial value problem (1.4) satisfies a similar growth estimate as
(2.1) in its existence time interval, so that it cannot blow-up in finite time.
Thus, the problem (1.4) is globally well-posed and the dynamical
system $S_1(t)$ is well-defined for all $t\in\mathbb{R}$. We thus have the
following result which is implicitly stated in \cite{Guo}:
\medskip

{\bf Theorem 2.4}\ \ {\it There exists a bounded and locally Lipschitz
continuous mapping $\Psi_1: H^{1/4}(\mathbb{R})\to C(\mathbb{R},H^{1/4}
(\mathbb{R}))$ such that for any $u_0\in H^{1/4}(\mathbb{R})$, the function
$u=\Psi_1(u_0)$ is a solution $($in distribution sense$)$ of the initial
value problem $(1.4)$, and it defines a dynamical system $S_1(t)$ in
$H^{1/4}(\mathbb{R})$ of the mKdV equation $(1.2)$. Moreover, there exists
a constant $C>0$ such that
$$
  \|u(\cdot,t)\|_{H^{1/4}(\mathbb{R})}\leq
  C(1+|t|)\|u_0\|_{H^{1/4}(\mathbb{R})} \quad
  \mbox{for all}\;\;t\in\mathbb{R}.
\eqno{(2.5)}
$$
$\Box$}

In Section 4 we shall use this theorem to prove Theorem 1.6.

The Miura transform is not a surjection from $H^{1/4}(\mathbb{R})$ to
$H^{-3/4}(\mathbb{R})$, cf. \cite{CCT}. Thus, we cannot use the relation
(2.4) and the weak continuity of $S_1(t)$ to get weak continuity of $S(t)$.
In order to prove Theorem 1.4 we shall use a generalized version of the Miura
transform --- the generalized Miura transform introduced by Christ,
Colliander and Tao \cite{CCT}, which is the mapping $(v,w)\mapsto u=M_g(v,w)$
defined by
$$
  u=\partial_x v+v^2+w.
\eqno{(2.6)}
$$
It can be easily verified that if $(v,w)$ is a solution of the initial value
problem
$$
\left\{
\begin{array}{l}
  \partial_t v+\partial_x^3 v=6(v^2+w)\partial_x v, \quad
  x\in\mathbb{R},\;\;\; t\in\mathbb{R},\\
  \partial_t w+\partial_x^3 w=6(v^2+w)\partial_x w, \quad
  x\in\mathbb{R},\;\;\; t\in\mathbb{R},\\
  v(x,0)=v_0(x), \quad x\in\mathbb{R}, \\
  w(x,0)=w_0(x), \quad x\in\mathbb{R},
\end{array}
\right.
\eqno{(2.7)}
$$
and $u_0=v_0'+v_0^2+w_0$, then $u=M_g(v,w)$ is a solution of the problem (1.3).
Using Lemma 2.2, we see immediately that $M_g$ maps $H^{1/4}(\mathbb{R})
\times H^1(\mathbb{R})$ into $H^{-3/4}(\mathbb{R})$, and it is a bounded
and locally Lipschitz continuous mapping. An important feature of $M_g$ is
that it is a surjection from  $H^{1/4}(\mathbb{R})\times H^1(\mathbb{R})$ to
$H^{-3/4}(\mathbb{R})$. More precisely, we have:
\medskip

{\bf Lemma 2.5}\ \ {\em For any $A>0$ there exists a Lipschitz continuous
mapping $W_A:H^{-3/4}(\mathbb{R})\to H^{1/4}(\mathbb{R})\times
H^1(\mathbb{R})$ such that $M\circ W_A=id$ when restricted to the ball
$B_A=\{u\in H^{-3/4}(\mathbb{R}):\|u\|_{H^{-3/4}}\leq A\}$.}
\medskip

{\em Proof}:\ \ See Lemma 10.1 of \cite{CCT}. $\quad\Box$
\medskip

In \cite{CCT} it was proved that the initial value problem (2.6) is
locally well-posed in $H^{1/4}(\mathbb{R})\times H^1(\mathbb{R})$. Since we
shall use this result later on, in the sequel we review some details of its
proof.

Let $\chi_k$ ($k=0,\pm 1,\pm 2,\cdots$) be the functions introduced in the
beginning of this section, and let
$$
  P_ku=\mathcal{F}^{-1}(\chi_k(\xi)\widehat{u}(\xi)) \quad
  \mbox{for}\;\; u\in\mathcal{S}'(\mathbb{R})
$$
($k=0,\pm 1,\pm 2,\cdots$). Given a time interval $I=[t_0,t_0+T]$, let $X=
X(\mathbb{R}\times I)$ be as before, and define the space $X^\ast=
X^\ast(\mathbb{R}\times I)$ by setting the norm
$$
  \|u\|_{X^\ast}=\|u\|_X+\Big(\sum_{k=-\infty}^\infty\|P_k u\|_X^2\Big)^{1/2}.
\eqno{(2.8)}
$$
Next we define the space $X^{\ast\ast}=X^{\ast\ast}(\mathbb{R}\times I)$ for
vector functions $(v,w)$ by setting the norm
$$
  \|(v,w)\|_{X^{\ast\ast}}=\|v\|_{X^\ast}+\|J_x^{3/4}w\|_{X^\ast}.
\eqno{(2.9)}
$$
Since $X\subseteq C(I,H^{1/4}(\mathbb{R}))$ and the embedding mapping is
continuous, we see easily that
$X^{\ast\ast}\subseteq C(I,H^{1/4}(\mathbb{R})\times H^1(\mathbb{R}))$, and
the embedding mapping is continuous. The local well-posedness result for the
problem (2.7) is as follows:
\medskip

{\bf Proposition 2.6}\ \ {\em Let $t_0=0$. For any $(v_0,w_0)\in H^{1/4}
(\mathbb{R})\times H^1(\mathbb{R})$ there exists $T=T(\|v_0\|_{H^{1/4}
(\mathbb{R})},\|w\|_{H^1(\mathbb{R})})>0$, such that the problem $(2.7)$ has
a unique solution $(v,w)$ in the space $X^{\ast\ast}$, and the mapping
$(v_0,w_0)\mapsto (v,w)$ from $H^{1/4}(\mathbb{R})\times H^1(\mathbb{R})$ to
$X^{\ast\ast}$ is locally Lipschitz continuous.}
\medskip

{\em Proof}:\ \ See Proposition 1 in \cite{CCT} and its proof. $\quad\Box$
\medskip

By Proposition 2.6, it follows that there exists a local dynamical system
$S^{\ast\ast}(t)$ in $H^{1/4}(\mathbb{R})\times H^1(\mathbb{R})$ for the
system of equations (2.7). Since the generalized Miura
transform $M_g$ maps a solution of (2.7) into a solution of (1.3) with $u_0
=v_0'+v_0^2+w_0$, by a similar argument as in the proof of Lemma 2.2 it
follows that
$$
  M_g[S^{\ast\ast}(t)(v_0,w_0)]=S(t)M_g(v_0,w_0)
\eqno{(2.10)}
$$
for any $(v_0,w_0)\in H^{1/4}(\mathbb{R})\times H^1(\mathbb{R})$ and any
$t\in\mathbb{R}$ such that $S^{\ast\ast}(t)(v_0,w_0)$ makes sense. In
Section 5 we shall use this relation and Proposition 2.6 to prove Theorem 1.5.

\section{Proof of Theorem 1.6}

\hskip 2em
In this section we give the proof of Theorem 1.6.

We denote $I=[-T,T]$. Let $u_{0n}\in H^{1/4}(\mathbb{R})$,
$n=1,2,\cdots$, and $u_0\in H^{1/4} (\mathbb{R})$ be such that
$u_{0n}\to u_0$ weakly in $H^{1/4}(\mathbb{R})$ as $n\to\infty$.
Then there exists constant $M>0$ such that
$$
  \|u_{0n}\|_{H^{1/4}}\leq M, \quad n=1,2,\cdots, \quad
  \mbox{and} \quad \|u_0\|_{H^{1/4}}\leq M.
\eqno{(3.1)}
$$
Let $u_n(x,t)=[S_1(t)u_{0n}](x)$, $n=1,2,\cdots$, and $u(x,t)=[S_1(t)u_0](x)$.
Let $T>0$ be given, and set $M_1=C(1+T)M$, where $C$ is the constant
appearing in (2.5). Then we have
$$
  \|u_n(\cdot,t)\|_{H^{1/4}}\leq M_1, \quad n=1,2,\cdots, \quad
  \mbox{and} \quad \|u(\cdot,t)\|_{H^{1/4}}\leq M_1
\eqno{(3.2)}
$$
for all $t\in I$. Using the equation (1.4), we further obtain
$$
  \|\partial_t u_n(\cdot,t)\|_{H^{-11/4}}\leq M_2, \quad n=1,2,\cdots, \quad
  \mbox{and} \quad \|\partial_t u(\cdot,t)\|_{H^{-11/4}}\leq M_2
\eqno{(3.3)}
$$
for all $t\in I$. Indeed, by (1.4) we have
$$
\begin{array}{rl}
  \|\partial_t u_n(\cdot,t)\|_{H^{-11/4}}\leq &
  \|\partial_x^3 u_n(\cdot,t)\|_{H^{-11/4}}+2
  \|\partial_x u_n^3(\cdot,t)\|_{H^{-11/4}}\\ [0.2cm]
\leq &
  \|u_n(\cdot,t)\|_{H^{1/4}}+C
  \|u_n^3(\cdot,t)\|_{H^{-5/6}}.
\end{array}
$$
Since $H^{1/4}(\mathbb{R})\subseteq L^4(\mathbb{R})$,
$L^{4/3}(\mathbb{R})\subseteq H^{-5/6}(\mathbb{R})$, and the
embeddings are continuous, we have
$$
  \|u_n^3(\cdot,t)\|_{H^{-5/6}}\leq C\|u_n^3(\cdot,t)\|_{L^{4/3}}
  =C\|u_n(\cdot,t)\|_{L^4}^3\leq C\|u_n(\cdot,t)\|_{H^{1/4}}^3.
$$
Hence
$$
  \|\partial_t u_n(\cdot,t)\|_{H^{-11/4}}\leq
  \|u_n(\cdot,t)\|_{H^{1/4}}+C\|u_n(\cdot,t)\|_{H^{1/4}}^3\leq M_2
$$
for all $t\in I$. The proof of the last inequality in (3.3) is
similar. We note that (3.2) also implies that for any $2\leq p\leq
4$,
$$
  \|u_n\|_{L^p(\mathbb{R}\times I)}\leq M_p, \quad
  n=1,2,\cdots.
\eqno{(3.4)}
$$
In addition, by the local well-posedness result for the problem
(1.4) that we reviewed in Section 2, from (3.2) we can also get the
following estimate
$$
  \|u_n\|_{X(\mathbb{R}\times I)}\leq C(T,M_1), \quad
  n=1,2,\cdots.
\eqno{(3.5)}
$$
We note that to get this estimate we need to divide the interval $I$
into small subintervals, with the number of them depending only on
$T$ and $M_1$.

By (3.2), (3.3) and a standard compactness result,  it follows that
there exists a subsequence $\{u_{n_k}\}$ of $\{u_n\}$ and a function
$u'\in L^2_{\rm loc}(\mathbb{R}\times I)$, such that for any $R>0$,
$$
  u_{n_k}\to u' \quad \mbox{strongly in}\;\; L^2([-R,R]\times I).
$$
This further implies, after passing to a subsequence when necessary,
that
$$
  u_{n_k}\to u' \quad \mbox{almost everywhere in}\;\;
  \mathbb{R}\times I.
\eqno{(3.6)}
$$
By (3.4), we have $u'\in L^p(\mathbb{R}\times I)$ for any $2\leq
p\leq 4$, and by (3.5), we also have $u'\in\widetilde{X}
(\mathbb{R}\times I)$, where $\widetilde{X}(\mathbb{R}\times I)$
denotes the function space of all measurable functions $u$ on
$\mathbb{R}\times I$ such that $\|u\|_{X(\mathbb{R}\times
I)}<\infty$. The last assertion is ensured by the fact that
$\widetilde{X}(\mathbb{R}\times I)$ is a $L^\infty$-type space,
i.e., it is the dual of a separable Banach space.

From (3.4) we see that
$$
  \{u_n\}\;\; \mbox{is bounded in}\;\; L^p(\mathbb{R}\times I),
  \;\;  2\leq p\leq 4,
\eqno{(3.7)}
$$
which further implies that
$$
  \{u_n^3\}\;\; \mbox{is bounded in}\;\; L^p(\mathbb{R}\times I),
  \;\; 1\leq p\leq \frac{4}{3}.
\eqno{(3.8)}
$$
Hence, by using the Vitali convergence theorem (see Corollary A.2 of
\cite{CK}), we infer from (3.6)--(3.8) that for any finite $R>0$,
$$
  u_{n_k}\to u '\quad \mbox{strongly in}\;\; L^p([-R,R]\times I),
  \;\;  1\leq p<4,
\eqno{(3.9)}
$$
$$
  u_{n_k}^3\to (u ')^3\quad \mbox{strongly in}\;\; L^p([-R,R]\times I),
  \;\;  1\leq p<\frac{4}{3}
\eqno{(3.10)}
$$
From (3.7)--(3.10) and the density of $C_0^\infty(\mathbb{R}\times
(-T,T))$ in $(L^p(\mathbb{R}\times I))'=L^{p'}(\mathbb{R}\times I)$
for $1<p<\infty$, we deduce that
$$
  u_{n_k}\to u' \quad \mbox{weakly in}\;\; L^p(\mathbb{R}\times I),
  \;\;  1<p<4,
$$
$$
  u_{n_k}^3\to (u')^3 \quad \mbox{weakly in}\;\; L^p(\mathbb{R}\times I),
  \;\;  1<p<\frac{4}{3}.
$$
Thus, by letting $k\to\infty$ in the equation
$$
  u_{n_k}(\cdot,t)= W(t)u_{0n_k}+2\partial_x\!
  \int_0^t  W(t-t')u_{n_k}^3(\cdot,t')dt', \quad
  k=1,2,\cdots,
$$
we see that $u'$ satisfies the equation
$$
  u'(\cdot,t)= W(t)u_0+2\partial_x\!
  \int_0^t  W(t-t')(u')^3(\cdot,t')dt'.
$$
Thus, since $u'\in\widetilde{X}(\mathbb{R}\times I)$ and, from the
proofs of Theorems 2.3 and 2.4 of \cite{KPV1} we see that the
solution of this equation in $\widetilde{X}(\mathbb{R}\times I)$ is
unique, we conclude that $u'=u$. Thus, we have proved that $\{u_n\}$
has a subsequence $\{u_{n_k}\}$ converging to $u$ almost everywhere
in $\mathbb{R}\times I$. Since the above argument works when the
sequence $\{u_n\}$ is replaced by any of its subsequences, it
follows that the following assertion holds:
$$
  u_n\to u \quad \mbox{almost everywhere in}\;\; \mathbb{R}\times I.
\eqno{(3.11)}
$$
As a consequence of this assertion, (4.10) (arbitrarily take $8/7<p\leq 4/3$) and the Vitali
convergence theorem, we see that also the following assertion holds:
$$
  u_n^3\to u^3\quad \mbox{strongly in}\;\; L^{8/7}([-R,R]\times I)\;\;
  \mbox{for any}\;\; R>0.
\eqno{(3.12)}
$$

Next, by (3.1), (3.2) and the density of $S(\mathbb{R})$ in
$H^{1/4}(\mathbb{R})$, we see that in order to prove Theorem 1.6 it
suffices to prove that (1.6) holds for any $\varphi\in
S(\mathbb{R})$. Since
$$
  (u,\varphi)_{H^{1/4}}=(u,\psi)_{L^2}, \quad
  \forall\varphi\in S(\mathbb{R}),
$$
where $\psi=\mathcal{F}_1^{-1}[(1+|\xi|^2)^{1/4}\mathcal{F}_1(\varphi)]\in
S(\mathbb{R})$, it follows that in order to prove that (1.6) holds for any
$\varphi\in S(\mathbb{R})$, it suffices to prove that the following holds
for any $\varphi\in S(\mathbb{R})$:
$$
  \lim_{n\to\infty}\sup_{|t|\leq T}|(u_n(\cdot,t)-u(\cdot,t),\varphi)_{L^2}|=0.
\eqno{(3.13)}
$$
Let
$$
  v_n(\cdot,t)= W(t)u_{n0}, \quad
  w_n(\cdot,t)=2\partial_x\!\int_0^t\! W(t-t')u_n^3(\cdot,t')dt', \quad
  n=1,2,\cdots,
$$
$$
  v(\cdot,t)= W(t)u_{0}, \quad
  w(\cdot,t)=2\partial_x\!\int_0^t\! W(t-t')u^3(\cdot,t')dt'.
$$
Then $u_n(\cdot,t)-u(\cdot,t)=[v_n(\cdot,t)-v(\cdot,t)]+[w_n(\cdot,t)-
w(\cdot,t)]$, $n=1,2,\cdots$. It can be easily shown that (cf. the proof of
Assertion 2 in Section 2.3 of \cite{CK})
$$
  \lim_{n\to\infty}\sup_{|t|\leq T}|(v_n(\cdot,t)-v(\cdot,t),\varphi)_{L^2}|=0
$$
  for any $\varphi\in S(\mathbb{R})$. Thus (3.13) follows if we prove that
$$
  \lim_{n\to\infty}\sup_{|t|\leq T}|(w_n(\cdot,t)-w(\cdot,t),\varphi)_{L^2}|=0
\eqno{(3.14)}
$$
for any $\varphi\in S(\mathbb{R})$. Let
$$
  z(\cdot,t)=\int_0^t\! W(t-t')u^3(\cdot,t')dt'
  \quad \mbox{and} \quad
  z_n(\cdot,t)=\int_0^t\! W(t-t')u_n^3(\cdot,t')dt'
$$
($n=1,2,\cdots$). Since
$w_n(\cdot,t)-w(\cdot,t)=2\partial_x[z_n(\cdot,t)- z(\cdot,t)]$,
$n=1,2,\cdots$, we see that (3.14) follows if we prove that
$$
  \lim_{n\to\infty}\sup_{|t|\leq T}|(z_n(\cdot,t)-z(\cdot,t),\varphi)_{L^2}|=0.
\eqno{(3.15)}
$$
The proof of this relation follows from a similar argument as in the proof of
Assertion 3 in Section 2.3 of \cite{CK}. Indeed, let $\chi_R$ be the
characteristic function of the interval $[-R,R]$ (for $x$ variable), and
denote
$$
  z^{(1)}_n(\cdot,t)=\!\int_0^t\! W(t-t')(\chi_R u_n^3(\cdot,t'))dt', \quad
  z^{(2)}_n(\cdot,t)=\!\int_0^t\!  W(t-t')((1-\chi_R)u_n^3(\cdot,t'))dt',
$$
$$
  z^{(1)}(\cdot,t)=\!\int_0^t\! W(t-t')(\chi_R u^3(\cdot,t'))dt', \quad
  z^{(2)}(\cdot,t)=\!\int_0^t\!  W(t-t')((1-\chi_R)u^3(\cdot,t'))dt'.
$$
Then we have
$$
\begin{array}{rl}
  (z_n(\cdot,t)&-z(\cdot,t),\varphi)_{L^2}=\displaystyle
  \int_{-\infty}^\infty[z_n(x,t)-z(x,t)]\varphi(x)dx
\\ [0.3cm]
  &=\displaystyle\int_{-\infty}^\infty[z^{(1)}_n(x,t)-z^{(1)}(x,t)]\varphi(x)dx
   +\int_{-\infty}^\infty[z^{(2)}_n(x,t)-z^{(2)}(x,t)]\varphi(x)dx\\ [0.3cm]
  &\equiv I_n^R(t)+J_n^R(t).
\end{array}
$$
Let
$$
  f_n^R(x,t)=\chi_R(x)u_n^3(x,t), \quad  f^R(x,t)=\chi_R(x)u^3(x,t).
$$
Then by the Cauchy inequality and the inhomogeneous Strichatz estimate we have
(note that $8/7$ is the dual of $8$ and $(2,\infty)$, $(8,8)$ are admissible
pairs)
$$
  \sup_{0\leq t\leq T}|I_n^R(t)|\leq
  \sup_{0\leq t\leq T}\|z^{(1)}_n(t)-z^{(1)}(t)\|_2\cdot\|\varphi\|_2\leq
  C\|f_n^R-f^R\|_{L^{8\over 7}(\mathbb{R}\times I)}\|\varphi\|_2
$$
($n=1,2,\cdots$). Hence, by (4.14) we have
$$
  \lim_{n\to\infty}\sup_{0\leq t\leq T}|I_n^R(t)|=0 \quad
  \mbox{for any fixed}\;\; R>0.
\eqno{(3.16)}
$$
Next, we compute
$$
\begin{array}{rl}
  J_n^R(t)&=\displaystyle\int_{-\infty}^\infty[z^{(2)}_n(x,t)-z^{(2)}(x,t)]
  \varphi(x)dx
\\ [0.3cm]
  &=\displaystyle\int_{-\infty}^\infty\int_0^t W(t-t')
   \{(1-\chi_R)[u_n^3(\cdot,t')-u^3(\cdot,t')]\}\varphi(x)dt' dx
\\ [0.3cm]
  &=\displaystyle\int_0^t \int_{-\infty}^\infty [1-\chi_R(x)]
  [u_n^3(x,t')-u^3(x,t')]\cdot W(t-t')\varphi
  dx dt'
\\ [0.3cm]
  &=\displaystyle\int_0^t \int_{-\infty}^\infty
  [u_n^3(x,t')-u^3(x,t')]\cdot
  [1-\chi_R(x)]W(t-t')\varphi dx dt'.
\end{array}
$$
From this expression and the H\"{o}lder inequality we have
$$
\begin{array}{rl}
  \displaystyle\sup_{0\leq t\leq T}|J_n^R(t)|
  \leq &\displaystyle\Big(\int_0^T\!\!\int_{-\infty}^\infty
  \Big|u_n^3(x,t')-u^3(x,t')\Big|^{8/7}
  dxdt'\Big)^{7/8}\Big(\int_0^T\!\!\int_{|x|\geq R}|W(t)\varphi|^8
  dx dt\Big)^{1/8}.
\end{array}
$$
By (4.10) (take $p=8/7$), the first term on the right-hand side is bounded by
a constant independent of $n$. Besides, since $(8,8)$ is an admissible pair, we
have $\|W(t)\varphi\|_{L^8(\mathbb{R}\times I)}\leq C\|\varphi\|_2$,
so that
$$
  \lim_{R\to\infty}\int_0^T\!\!\int_{|x|\geq R}|W(t)\varphi|^8
  dx dt=\lim_{R\to\infty}\int_{|x|\geq R}\int_0^T|W(t)\varphi|^8
  dx dt=0.
$$
Hence
$$
  \lim_{R\to\infty}\sup_{0\leq t\leq T}|J_n^R(t)|=0 \quad
  \mbox{uniformly for}\;\; n\in\mathbb{N}.
\eqno{(3.17)}
$$
By (3.16) and (3.17), we obtain (3.15). This completes the proof of
Theorem 1.6. $\quad\Box$

\section{Proof of Theorem 1.5}

\hskip 2em
In this section we give the proof of Theorem 1.5.

We first prove a preliminary lemma. Let $M_g$ be the generalized Miura
transform given by (2.6) and $W_A$ be the Lipschitz continuous mapping
from the ball $B_A=\{\phi\in H^{-3/4}(\mathbb{R}):
\|\phi\|_{H^{-3/4}}\leq A\}$ of $H^{-3/4}(\mathbb{R})$ to $H^{1/4}(\mathbb{R})\times H^1(\mathbb{R})$
ensured by Lemma 2.5 such that
$$
  M_g\circ W_A=id.
$$
\medskip

{\bf Lemma 4.1}\ \ {\em $(i)$ Let $(v_n,w_n)\in
H^{1/4}(\mathbb{R})\times H^1(\mathbb{R})$, $n=1,2,\cdots$, and
$(v,w)\in H^{1/4}(\mathbb{R})\times H^1(\mathbb{R})$. Let
$u_n=M_g(v_n,w_n)$, $n=1,2,\cdots$, and $u=M_g(v,w)$. Assume that
$(v_n,w_n)\to (v,w)$ weakly in $H^{1/4}(\mathbb{R})\times
H^1(\mathbb{R})$. Then $u_n\to u$ weakly in $H^{-3/4}(\mathbb{R})$.

$(ii)$ Conversely, let $u_n\in H^{-3/4}(\mathbb{R})$, $n=1,2,\cdots$, and $u
\in H^{-3/4}(\mathbb{R})$. Let $(v_n,w_n)=W_A(u_n)$, $n=1,2,\cdots$, and
$(v,w)=W_A(u)$. Assume that $u_n\to u$ weakly in $H^{-3/4}(\mathbb{R})$. Then
$(v_n,w_n)\to (v,w)$ weakly in $H^{1/4}(\mathbb{R})\times H^1(\mathbb{R})$.}
\medskip

{\em Proof}:\ \ $(i)$ We write
$$
  u_n=\partial_x v_n+v_n^2+w_n \quad \mbox{and} \quad u=\partial_x v+v^2+w.
$$
Clearly, $\partial_x v_n\to\partial_x v$ weakly in $H^{-3/4}(\mathbb{R})$
and $w_n\to w$ weakly in $H^{-3/4}(\mathbb{R})$. Thus, we only need to prove
that $v_n^2\to v^2$ weakly in $H^{-3/4}(\mathbb{R})$, which is almost obvious.
Indeed, by the compact embedding $H^{1/4}(\mathbb{R})\hookrightarrow L^2[a,b]$
for any $-\infty<a<b<\infty$, we easily see that $v_n^2\to v^2$ strongly in
$L^1[a,b]$ for any $-\infty<a<b<\infty$, so that for any $\varphi\in
C_0^\infty(\mathbb{R})$ we have
$$
  |(v_n^2-v^2,\varphi)_{L^2}|\leq\|v_n^2-v^2\|_{L^1[a,b]}
  \|\varphi\|_\infty\to 0 \;\;\;\mbox{as}\;\; n\to\infty,
$$
where $a,b$ are real numbers such that $\supp\varphi\subseteq [a,b]$. This
implies that for any $\psi\in J_x^{3/2}(C_0^\infty(\mathbb{R}))\subseteq
H^{-3/4}(\mathbb{R})$ we have
$$
  (v_n^2-v^2,\psi)_{H^{-3/4}}\to 0 \;\;\;\mbox{as}\;\; n\to\infty.
$$
Since $C_0^\infty(\mathbb{R})$ is dense in $H^{3/4}(\mathbb{R})$ and
$J_x^{3/2}$ is an isomorphism of $H^{3/4}(\mathbb{R})$ onto
$H^{-3/4}(\mathbb{R})$, we see that $J_x^{3/2}(C_0^\infty(\mathbb{R}))$ is
dense in $H^{-3/4}(\mathbb{R})$. Thus, by the boundedness of the sequence
$\{v_n^2-v^2\}$ in $H^{-3/4}(\mathbb{R})$, we conclude that the above relation
holds for all $\psi\in H^{-3/4}(\mathbb{R})$. This proves the assertion $(i)$.

$(ii)$ We first recall the construction of the mapping $W_A:H^{-3/4}(\mathbb{R}
\to H^{1/4}(\mathbb{R})\times H^1(\mathbb{R})$ (cf. the proof of Lemma 10.1
of \cite{CCT}). We know that there are sufficiently large constants $C_A,C_A'
>0$ such that, denoting by $P$ the Fourier projection to the frequency region
$|\xi|\geq C_A$, for any given $\phi\in H^{-3/4}(\mathbb{R})$ such that
$\|\phi\|_{H^{-3/4}}\leq A$, the mapping $\varphi\mapsto\partial_x^{-1}P(\phi-
\varphi^2)$ is a contraction on the ball $\|\varphi\|_{H^{1/4}}\leq C_A'A$ in
$H^{1/4}(\mathbb{R})$. Let $\varphi=L(\phi)$ be the unique fixed point of
this mapping and $\psi=(I-P)(\phi-\varphi^2)$. Then we have $W_A(\phi)=
(\varphi,\psi)$.

We now proceed to prove the assertion $(ii)$. Since $u_n\to u$ weakly in
$H^{-3/4}(\mathbb{R})$, $\{u_n\}$ is bounded in $H^{-3/4}(\mathbb{R})$, which
implies that $\{v_n\}$ is bounded in $H^{1/4}(\mathbb{R})$. It follows that
there exists subsequence of $\{v_n\}$ which weakly converges in $H^{1/4}
(\mathbb{R})$. Thus, to prove that $v_n\to v$ weakly in $H^{1/4}(\mathbb{R})$
we only need to prove that if a subsequence
of $\{v_n\}$ weakly converges to an element $v'\in H^{1/4}(\mathbb{R})$, then
$v'=v$. For simplicity of notation we assume that the whole sequence
$v_n\to v'$ weakly in $H^{1/4}(\mathbb{R})$. From the proof of $(i)$ we see
that this implies that $v_n^2\to (v')^2$ weakly in $H^{-3/4}(\mathbb{R})$.
For every $n\in\mathbb{N}$ we have
$$
  v_n=\partial_x^{-1}P(u_n-v_n^2),
$$
or equivalently,
$$
  \partial_x v_n=P(u_n)-P(v_n^2).
$$
Letting $n\to\infty$ and considering the weak limits, we get
$$
  \partial_x v'=P(u)-P((v')^2).
$$
This shows that $v'$ is a fixed point of the mapping $\varphi\mapsto
\partial_x^{-1}P(u-\varphi^2)$. Since $W_A(u)=(v,w)$, we see that $v$ is also
a fixed point of this mapping. By uniqueness of the fixed point, we obtain
$v'=v$. Hence, the desired assertion follows. From this assertion and the
relations $w_n=(I-P)(u_n-v_n^2)$ ($n=1,2,\cdots$) and $w=(I-P)(u-v^2)$ it
follows immediately that $w_n\to w$ weakly in $H^1(\mathbb{R})$. This
completes the proof. $\quad\Box$
\medskip

We are now ready to prove Theorem 1.5.

Let $u_{0n}\in H^{-3/4}(\mathbb{R})$, $n=1,2,\cdots$, and $u_0\in H^{-3/4}
(\mathbb{R})$ be such that $u_{0n}\to u_0$ weakly in $H^{-3/4}(\mathbb{R})$
as $n\to\infty$. Then there exists constant $M>0$ such that
$$
  \|u_{0n}\|_{H^{-3/4}}\leq M, \quad n=1,2,\cdots, \quad
  \mbox{and} \quad \|u_0\|_{H^{-3/4}}\leq M.
\eqno{(4.1)}
$$
Let $u_n(x,t)=[S(t)u_{0n}](x)$, $n=1,2,\cdots$, and $u(x,t)=[S(t)u_0](x)$.
Given $T>0$ be given, we set $A=C(1+T)M$, where $C$ is the constant appearing
in (2.1). Then we have
$$
  \|u_n(\cdot,t)\|_{H^{-3/4}}\leq A, \quad n=1,2,\cdots, \quad
  \mbox{and} \quad \|u(\cdot,t)\|_{H^{-3/4}}\leq A
\eqno{(4.2)}
$$
for all $t\in I$. Let $(v_{0n},w_{0n})=W_A(u_{0n})$, $n=1,2,\cdots$,
and $(v_0,w_0)=W_A(u_0)$. Let
$$
  M=M_A:=\sup_{(\varphi,\psi)\in W_A(B_A)}\max\{\|\varphi\|_{H^{1/4}},
  \|\psi\|_{H^1}\}.
$$
By Proposition 2.6, for this constant $M$ there is a corresponding constant
$\delta>0$ such that for any $t_0\in\mathbb{R}$ and any $(\varphi,\psi)\in
H^{1/4}(\mathbb{R})\times H^1(\mathbb{R})$ satisfying $\|\varphi\|_{H^{1/4}}
\leq M$ and $\|\psi\|_{H^1}\leq M$, the initial value problem
$$
\left\{
\begin{array}{l}
  \partial_t v+\partial_x^3 v=6(v^2+w)\partial_x v, \quad
  x\in\mathbb{R},\;\;\; t_0-\delta\leq t\leq t_0+\delta,\\
  \partial_t w+\partial_x^3 w=6(v^2+w)\partial_x w, \quad
  x\in\mathbb{R},\;\;\; t_0-\delta\leq t\leq t_0+\delta,\\
  v(x,t_0)=\varphi(x), \quad x\in\mathbb{R}, \\
  w(x,t_0)=\psi(x), \quad x\in\mathbb{R},
\end{array}
\right. \eqno{(4.3)}
$$
has a unique solution $(v,w)$ in the space $X^{\ast\ast}=X^{\ast\ast}
(\mathbb{R}\times [t_0-\delta,t_0+\delta])$ (see (2.9) for the definition of
this space), which we also denote as $(v(\cdot,t),w(\cdot,t))=
S^{\ast\ast}_{t_0}(t)(\varphi,\psi)$. Moreover, we have the estimate
$$
  \|(v,w)\|_{X^{\ast\ast}}\leq C(\max\{\|\varphi\|_{H^{1/4}},
  \|\psi\|_{H^1}\})\leq C(M),
\eqno{(4.4)}
$$
where $C:[0,\infty)\to [0,\infty)$ is a nondecreasing function. Let
$$
  (v_n(\cdot,t),w_n(\cdot,t))=S^{\ast\ast}_0(v_{0n},w_{0n}),
  \quad n=1,2,\cdots,
$$
and $(v(\cdot,t),w(\cdot,t))=S^{\ast\ast}_0(v_0,w_0)$. Since $(v_{0n},w_{0n})
\in W_A(B_A)$, we have $\|v_{0n}\|_{H^{1/4}}\leq M$ and $\|w_{0n}\|_{H^1}\leq
M$, $n=1,2,\cdots$. Thus
$$
  \|(v_n,w_n)\|_{X^{\ast\ast}}\leq C(M), \quad n=1,2,\cdots.
\eqno{(4.5)}
$$
Moreover, since $u_{0n}\to u_0$ weakly in $H^{-3/4}(\mathbb{R})$ as
$n\to \infty$, by Lemma 4.1 $(ii)$ we see that
$$
  (v_{0n},w_{0n})\to (v_0,w_0) \;\; \mbox{weakly in}\;\; H^{1/4}(\mathbb{R})
  \times H^1(\mathbb{R})\;\; \mbox{as}\;\; n\to\infty.
\eqno{(4.6)}
$$
From (4.5) and the definition of the space $X^{\ast\ast}$ it follows
immediately that
$$
  \|v_n\|_{L_T^\infty H_x^{1/4}}\leq C(T,M)
\eqno{(4.7)}
$$
and
$$
  \|w_n\|_{L_T^\infty H_x^1}\leq C(T,M)
\eqno{(4.8)}
$$
for all $n=1,2,\cdots$.

Using the estimates (4.8) and a similar argument as in the proof of
Theorem 1.6, we conclude that for any $\varphi\in
H^{1/4}(\mathbb{R})$ and $\psi\in H^1(\mathbb{R})$ we have
$$
  \lim_{n\to\infty}\{\sup_{|t|\leq\delta}
  |(v_n(\cdot,t)-v(\cdot,t),\varphi)_{H^{1/4}}|
  +\sup_{|t|\leq\delta}|(w_n(\cdot,t)-w(\cdot,t),\psi)_{H^1}|\}=0
$$
Using this relation, the relation (2.10) and a similar argument as
in the proof of Lemma 4.1 $(i)$, we obtain that for any $\phi\in
H^{-3/4}(\mathbb{R})$,
$$
  \lim_{n\to\infty}\sup_{|t|\leq\delta}
  |(u_n(\cdot,t)-u(\cdot,t),\phi)_{H^{-3/4}}|=0.
$$
Now let $m=T/\delta$ if $T/\delta$ is an integer, and
$m=[T/\delta]+1$ otherwise. We divide the time interval $I$ into
$2m$ subintervals $I_{\pm 1}$, $I_{\pm 2}$, $\cdots$, $I_{\pm m}$,
where
$$
  I_j=[(j\!-\!1)\delta,j\delta], \quad j=1,2,\cdots, m\!-\!1, \quad
  I_m=[(m\!-\!1)\delta,T],
$$
and $I_{-j}=-I_j$, $j=1,2,\cdots,m$. By inductively using the result we have
just proved to every pair of intervals $I_{\pm j}$, $j=1,2,\cdots,m$, we see
that for any $\phi\in H^{-3/4}(\mathbb{R})$,
$$
  \lim_{n\to\infty}\sup_{t\in I_j}
  |(u_n(\cdot,t)-u(\cdot,t),\phi)_{H^{-3/4}}|=0, \quad
  j=\pm 1,\pm 2,\cdots, \pm m.
$$
Hence, for any $\phi\in H^{-3/4}(\mathbb{R})$,
$$
  \lim_{n\to\infty}\sup_{|t|\leq T}
  |(u_n(\cdot,t)-u(\cdot,t),\phi)_{H^{-3/4}}|=0.
$$
This completes the proof of Theorem 1.5. $\quad\Box$

\end{document}